\documentclass{llncs}
\usepackage{amsmath}
\usepackage{amssymb}
\usepackage{eulervm}
\usepackage{enumerate}
\usepackage{varioref}
\usepackage{charter}

\title{Graphs without even holes or diamonds}
\author{
 Ton~Kloks%
\thanks{During this research this coauthor was 
a guest of the School of Computing, University of Leeds, 
Leeds LS2 9JT, UK.}  
}

\begin{document}

\maketitle

\begin{abstract}
An even hole is an induced chordless cycle of even length 
at least four. A diamond is an induced subgraph isomorphic to 
$K_4-e$. We show that graphs without even holes and without diamonds 
can be decomposed via clique-separators into graphs that have uniformly 
bounded cliquewidth. 
\end{abstract}

\section{Introduction}
%%%%%%%%%%%%%%%%%%%%%%

We consider undirected graphs without loops or multiple edges. 

\begin{definition}
A {\em chord\/} in a cycle is an edge connecting two vertices of the cycle 
that are not adjacent in the cycle. A cycle is {\em chordless\/} if 
it has no chord. 
A {\em hole\/} is a chordless cycle 
of length at least four. A graph is {\em chordal\/} if it has no holes.  
\end{definition}

Similarly, a {\em chordless path\/} in a graph $G$ 
is a set of vertices $P$ such that $G[P]$ is a path.  
A hole is {\em even\/} if it has even length. 
A {\em diamond\/} is an induced $K_4-e$. In this paper we consider 
graphs without even holes and diamonds. 

Graphs without even holes were studied 
in~\cite{kn:addario,kn:chudnovsky,kn:chudnovsky2,%
kn:conforti,kn:conforti2,kn:murilo,kn:da}. 
These graphs can be recognized in 
polynomial time~\cite{kn:chudnovsky,kn:conforti2}. 
It was shown that every graph without even holes has a 
{\em bisimplicial extreme\/}, that is a vertex whose neighborhood is the 
union of two cliques~\cite{kn:chudnovsky2}. 
Graphs without diamonds nor even holes were first 
studied in~\cite{kn:kloks2}. 
We showed in~\cite{kn:kloks2} that 
every graph without diamonds and without even holes has a 
{\em simplicial extreme\/}, that is a vertex which is either simplicial 
or which has degree 2.  

Cliquewidth was introduced by Courcelle and Olariu in~\cite{kn:courcelle}. 
For integer $k$ we define the following
{\em $k$-composers\/}:
\begin{enumerate}[\rm 1.]
\item For $i \in \{1,\dots,k\}$,
$i(x)$ creates a vertex $x$ with label $i$.
\item For distinct $i,j \in \{1,\dots,k\}$, let $\eta_{i,j}$
be the operation which adds all edges joining vertices with
label $i$ to vertices with label $j$.
\item For $i,j \in \{1,\dots,k\}$, let $\rho_{i \rightarrow j}$
be the operation which relabels every vertex with label $i$ with the  
label $j$.
\item The operation $G \oplus H$ creates the graph which is
the disjoint union
of two labeled graphs $G$ and $H$.
\end{enumerate}

\begin{definition}
A graph $G$ has {\em cliquewidth $k$\/} if $G$ can be constructed  
via a series of $k$-composers.
\end{definition}

The corresponding graph-decomposition is called a $k$-expression.  
Courcelle showed that those problems that can be 
formulated in monadic second-order logic without quantification over  
subsets of edges,   
can be solved efficiently for graphs of bounded 
cliquewidth (see~\cite{kn:courcelle3}). 
For these algorithms a $k$-expression is a necessary ingredient. 
However, this obstacle was taken away with the introduction of 
rankwidth. The 
graph classes for which 
the parameters 
cliquewidth and rankwidth are bounded coincide: 
It was shown in~\cite{kn:oum4,kn:oum2} that $rw \leq cw \leq 2^{rw+1}-1$. 
An algorithm that computes a rank-decomposition-tree  
for a graph of bounded rankwidth was 
obtained in
\cite{kn:hlineny2,kn:hlineny}. This algorithm runs in $O(n^3)$ time.  
Interestingly, it is still unknown whether 
{\sc cliquewidth} is fixed-parameter tractable 
(see, {\em e.g.\/},~\cite{kn:flum}). 
Computing cliquewidth is 
NP-complete~\cite{kn:fellows}. The NP-completeness of {\sc rankwidth}  
seems to follow from arguments given in~\cite{kn:hlineny}. 

In this paper we show that graphs without even holes and without 
diamonds can be decomposed via clique separators 
into graphs which have a uniform bound on the 
cliquewidth. 
It easily follows that graphs in this class can be 
recognized in $O(n^3)$ time. 
This result was anticipated by~\cite{kn:da}. That paper 
shows that planar graphs without even holes have uniformly 
bounded treewidth. 

\section{Birdcages and links}
%%%%%%%%%%%%%%%%%%%%%%%%%%%%%
\label{section birdcages}

For two sets $A$ and
$B$ we write $A+B$ for $A \cup B$ and $A-B$ for $A \setminus B$.
For a set $A$ and an element $x$ we also write $A+x$ instead of
$A+\{x\}$.
For a vertex $x$ we write $N(x)$ for the set of its neighbors and we
write $N[x]=x+N(x)$ for its {\em closed\/} neighborhood.
For a subset $S$ of vertices we write
$N(S)=\bigcup_{x \in S}N(x)-S$ and we write 
$N[S]=N(S)+S$. For a subset $S$ of vertices of
a graph $G$ we write $G[S]$ for the subgraph of $G$ {\em induced\/}
by $S$. For a graph $G=(V,E)$ and a subset $S$ of its vertices we write
$G-S$ for the graph $G[V-S]$. If $S$ consists of a single vertex
$x$ we also write $G-x$ instead of $G-\{x\}$.

Let $\mathcal{G}$ be the class of graphs without even holes
and without diamonds.

\begin{definition}
A {\em simplicial\/} is a vertex whose neighborhood 
induces a clique. A {\em simplicial extreme\/} is a vertex which 
either 
is simplicial, or has degree two. 
\end{definition}

\begin{lemma}[\cite{kn:kloks2}]
\label{simplicial extreme}
Every graph which has no even holes nor diamonds is either a 
clique or has two simplicial extremes that are not adjacent. 
\end{lemma}

\begin{definition}
Let $C_1$ and $C_2$ be cliques. Consider all chordless 
paths from a vertex of $C_1$ to a vertex of $C_2$ which does 
not use any other vertex of $C_1$ or $C_2$. If the lengths of 
all those 
paths have the same parity, then this collection of 
chordless paths is called a 
{\em $C_1,C_2$-link\/}.  
\end{definition}

In this section we analyze the structure of links. Our basic 
building blocks are birdcages. 

\begin{definition}
A {\em birdcage\/} is a graph with a specified  
clique $F$,  
called the {\em floor\/}, and a specified 
vertex $h \not\in F$, 
called the {\em hook\/}, and a collection of 
paths $\mathcal{P}$, one from each vertex in $F$ 
to $h$, such that
\begin{enumerate}[\rm (a)]
\item all path of $\mathcal{P}$ have odd length or all paths 
in $\mathcal{P}$ have even length, and  
\item either at most one path in $\mathcal{P}$ has length one 
or all paths in $\mathcal{P}$ have length one. 
\end{enumerate}
\end{definition}

We consider two ways to connect two birdcages. 

\begin{definition}
Let $B_i=(F_i,h_i)$, $i=1,2$ be two birdcages. A {\em join\/} 
of $B_1$ and $B_2$ is the graph obtained either by 
adding all edges between the floor or the hook 
of $B_1$ and the floor or the hook of $B_2$, or by identifying 
the two hooks $h_1$ and $h_2$. 
In case $|F_1|=|F_2|$ and if exactly one of $B_1$ and $B_2$ 
is a clique, then we also call the vertex-by-vertex  
identification of $F_1$ 
and $F_2$ a join. In case $|F_1|=|F_2|=2$ and if for 
$i=1,2$ there 
exists exactly one vertex $x_i \in F_i$ which is adjacent 
to $h_i$, then we also call the identification of 
$F_1$ and $F_2$ such that $x_1$ and $x_2$ are identified, 
a join.    
\end{definition}

\begin{definition}
Let $B=(F,h)$ be a birdcage. Let $h^{\ast}$ be a fixed  
neighbor of $h$ such that $h$ and $h^{\ast}$ have no 
neighbors in common. The edge $(h,h^{\ast})$ is called 
the {\em skew-edge\/} of $B$. 
\end{definition}

\begin{definition}
Let $B_i=(F_i,h_i)$, $i=1,2$ be two birdcages each with 
a skew-edge $(h_i,h_i^{\ast})$. A {\em skew-join\/} of 
$B_1$ and $B_2$ is the graph obtained by identifying 
each of $h_i$ and $h_i^{\ast}$ with one of $h_{3-i}$ and 
$h_{3-i}^{\ast}$. Let $B_1=(F_1,h_1)$ be a birdcage and let 
$B_2=(F_2,h_2)$ be a birdcage with a skew-edge $(h_2,h_2^{\ast})$. 
The identification of $h_1$ with one of $h_2,h_2^{\ast}$, and also 
the join of $h_1$ or $F_1$ with one of $h_2$ and $h_2^{\ast}$,  
is    
called a skew-join.  
\end{definition}

We consider one more operation. 

\begin{definition}
Let $B=(F,h)$ be a birdcage and let $L$ be a link from a 
vertex to a clique. Assume that the length of 
$L$ has the same parity as the lengths of the paths connecting 
$h$ with $F$ in $B$. A {\em replacement\/} is the operation of 
substituting one $h,F$-path in $B$ by $L$. 
The clique 
of $L$ is joined to the other vertices of $F$.  
\end{definition}

\noindent
We call a birdcage in which some paths have been replaced by links 
again a birdcage. 

\begin{theorem}
\label{link}
There exists a natural number $t$ such that every link has 
cliquewidth at most $t$. 
\end{theorem}
\begin{proof}
Consider a link $L$ between two nonadjacent vertices 
$x$ and $y$.  Assume first that $x$ has at least two nonadjacent 
neighbors. Between any pair of cliques $C_1$ and $C_2$ in 
$N(x)$ all chordless paths must be odd thus the induced 
$C_1,C_2$-paths in $L$ form a $C_1,C_2$-link. Note that $y$ is 
not in any of these induced sublinks since that implies an even 
hole. Let $L^{\ast}$ be the subgraph of $L$ induced by 
these sublinks. We prove that $L^{\ast}$ induces a birdcage with hook 
$x$ and some separating floor $F$. The floor $F$ separates $L^{\ast}$ 
from a link from $F$ to $y$. 
Consider the component $C$ of $L-L^{\ast}$ 
that contains $y$. We first show that $N(C)$ is a clique. 
Let $a$ and $b$ be nonadjacent vertices in $N(C)$. Since an $a,b$-path 
with internal vertices in $C$ cannot connect two vertices that are in a 
link, $a$ and $b$ must be two nonadjacent vertices in a birdcage. 
However, by induction 
this implies that there is a link (either between two cliques 
in $N(x)$ or from a clique in $N(x)$ to $y$) which contains a birdcage with 
two ends of the same type $\in \{\mbox{hook},\mbox{floor}\}$, which 
is a contradiction.  
The chordless paths from $y$ to any of the cliques in $N(x)$ 
induce a link. By induction, either there exists a join 
between two cliques $N(C)$ and $N(L^{\ast})$ and two links 
between $N(C)$ and $x$ and $N(L^{\ast})$ and $y$ or there exists 
a skew-edge $(a,b)$. The skew-edge connects three links; two 
to $N(x)$ and one to $y$. Note however that such a skew-join 
is only possible in case $N(x)$ is a clique (see below) since otherwise 
there would exist an even hole.    

In case $N(x)$ is a clique, we obtain a similar decomposition, except that 
in this case $y$ can be a vertex of $L^{\ast}$, which in this case 
is defined as follows. When $N(x)$ is a clique we 
consider the subgraph $L^{\ast}$ induced by chordless paths between 
different vertices of $N(x)$, 
with internal vertices in $L-N[x]$. The structure is proved by induction 
on $|N(x)|$.  
Note that a skew-connection  
between $L^{\ast}$ and $L-L^{\ast}$  
is only possible in case $N(x)$ is a clique.  
  
The induced subgraph $L^{\ast}$ is a {\em bag\/}. It follows that 
$L$ can be decomposed into a tree of bags glued together along 
clique cutsets either via joins or via skew-joins.   

Consider the {\em cross-edges\/} in an $x,y$-link $L$ for 
nonadjacent vertices $x$ and $y$: Assume two vertices $a$ and $b$ are 
adjacent but the edge is not an edge of $L$. 
Then $a$ and $b$ 
are contained in a birdcage, since one cannot be on a 
chordless $x,y$-path of the other. It follows that there exists a 
birdcage $B=(F,h)$ with a vertex $a \in F$ adjacent to $h$, and the 
vertex $a$ is adjacent to some cliques in the other (replaced) 
$h,F$-paths. 
By induction on the induced link-structures, the chordless paths 
from $a$ to either 
$x$ or $y$ must form a link. This implies that 
if $a$ is adjacent to two 
vertices in different paths of a birdcage 
$B^{\prime}=(F^{\prime},h^{\prime})$ then $a$ is also adjacent to 
$h^{\prime}$. Note also that for each birdcage $B^{\ast}$ 
there is at most 
one vertex $a$ (adjacent to the hook of $B^{\ast}$) which has 
crossing edges to paths of $B^{\ast}$.  
This proves that these birdcages with their    
cross-edges  
can be described by a bounded cliquewidth expression. 

By the recursive definition, this proves the theorem. 
\qed\end{proof}

\section{Link-extensions}
%%%%%%%%%%%%%%%%%%%%%%%%%

In this section we extend links in a recursive way. Consider two 
nonadjacent vertices $a$ and $b$ in an  
$x,y$-link $L$ induced by the collection of chordless $x,y$-paths 
in a graph $G \in \mathcal{G}$. Assume there exists a chordless 
$a,b$-path with internal vertices $\not\in L$.  
Then $a$ and 
$b$ must be vertices of a birdcage, since they cannot both lie on a 
common chordless $x,y$-path. Thus there exists a birdcage $B=(F,h)$ 
with a vertex $a \in F$ which is adjacent to $h$. There are 
(extended) links, ``cross-links,'' 
from $a$ to cliques or skew-joins in the other 
replaced paths in $B$. 
These cross-links act as the cross-edges introduced in the proof 
of Theorem~\ref{link}, except that 
the cross-links can also go to skew-edges. 
Note that if $a$ has links to 
cliques or skew-joins in different paths of some birdcage $B^{\prime}$ 
then $a$ must also be adjacent to the hook of $B^{\prime}$. Furthermore, 
for each constituent birdcage $B^{\ast}$ in $L$ 
there is at most one vertex $a$ 
with links to paths of $B^{\ast}$.        

\begin{theorem}
\label{link-extension}
There exists a number $t$ such that 
every graph $G \in \mathcal{G}$ without clique-separators 
has cliquewidth $t$. 
\end{theorem}
\begin{proof}
Assume $G$ is not a clique. We proved in~\cite{kn:kloks2} 
that $G$ has a vertex $\omega$ with two nonadjacent neighbors 
$x$ and $y$. The collection of chordless $x,y$-paths forms 
an $x,y$-link in $G-\omega$. Let $L^{\ast}$ be the extension of 
$L$ in $G-\omega$. 
By definition of a link-extension, and since there are no 
clique-separators, $G-\omega=L^{\ast}$. 
\qed\end{proof}

\section{Cliquewidth of (even-hole,diamond)-free graphs}
%%%%%%%%%%%%%%%%%%%%%%%%%%%%%%%%%%%%%%%%%%%%%%%%%%%%%%%%
  
\begin{definition}
A {\em splitgraph\/} is a graph $H=H(C,I,E)$ 
with a partition of the 
vertices in a clique $C$ and an independent set $I$.  
\end{definition}

\begin{lemma}[\cite{kn:boliac,kn:makowski}]
\label{unbounded split}
Cliquewidth is unbounded for splitgraphs. 
\end{lemma}
 
\begin{definition}
A {\em birdcage-split\/} is a graph which can be constructed from a 
splitgraph $H=H(C,I,E)$  
by replacing all edges incident with every 
vertex $x \in I$ either by an even or an odd $>1$-length path. 
\end{definition}

\noindent
Thus a birdcage-split $H$ consists of a clique $C$ and a 
collection $\mathcal{F}$ of subsets of $C$ and an independent set 
$\mathcal{H}$ of hooks. Each hook forms a birdcage in 
$H$ with floor $F \in \mathcal{F}$. 

Let $\mathcal{BS}$ be the class of birdcage-splits. 

\begin{lemma}
Birdcage-splits do not have even holes or diamonds. 
\end{lemma}
\begin{proof}
Every hole is contained in a birdcage, thus it is odd. 
Since all paths in birdcages have length more than one, there 
is no diamond. 
\qed\end{proof}

\begin{lemma}
\label{unbounded}
For every natural number $t$ there exists a graph in $\mathcal{BS}$ 
with cliquewidth more than $t$. 
\end{lemma}
\begin{proof}
Let $G=G(C,I,E)$ be a splitgraph and let $G^{\ast}$ be 
the birdcage-split 
obtained from $G$ by subdividing every edge $(x,y)$ 
with $x \in C$ and $y \in I$ by a single vertex. Let $z$ be a 
subdivision vertex of $G^{\ast}$ and consider a 
{\em local complementation\/} at $z$. If $x \in C$ and $y \in I$ are the 
two neighbors of $z$ then this adds the edge $(x,y)$ to $G^{\ast}$. 
Let $\Hat{G}$ be the graph obtained from $G^{\ast}$ by doing a local 
complementation at every subdivision vertex. Then $G$ is an induced 
subgraph of $\Hat{G}$, that is, $G$ is a {\em vertex-minor\/} 
of $G^{\ast}$. This implies that the rankwidth of $G^{\ast}$ is at least 
the rankwidth of $G$~\cite{kn:oum3}. By Lemma~\ref{unbounded split} this 
proves the claim.    
\qed\end{proof}

\begin{theorem}
There exists a natural number $t$ such that every $G \in \mathcal{G}$ 
either has a clique-separator or has cliquewidth at most $t$. 
\end{theorem}
\begin{proof}
Links and link-extensions 
can be generated by a recursive function of bounded width. 
The graph can be decomposed by clique-separators into these 
link-extensions. 
\qed\end{proof}

Recall that a decomposition by clique separators can be 
obtained in $O(n^3)$ time~\cite{kn:tarjan}. 

\begin{corollary}
Graphs in $\mathcal{G}$ can be recognized in $O(n^3)$ time. 
\end{corollary}

\section{Remarks on geodetic graphs}
%%%%%%%%%%%%%%%%%%%%%%%%%%%%%%%%%%%%

Geodetic graphs were introduced by Ore.

\begin{definition}[\cite{kn:bosak,kn:ore}]
A graph is {\em geodetic\/} if for every pair of vertices the
shortest path between them is unique.
\end{definition}

Note that a graph is geodetic if and only if for every vertex 
$x$, every vertex $y \in N_k(x)$\footnote{As usual,  
$N_k(x)=\{y\;|\; d(x,y)=k\}$.} is adjacent to exactly 
one vertex in $N_{k-1}(x)$, $k \geq 2$; see~\cite{kn:parthasarathy}. 
This settles 
the recognition problem. 
It follows from the definition that 
geodetic graphs have no induced 
diamond and no induced $C_4$. 
 
One partial characterization of geodetic graphs of
diameter two appeared in~\cite{kn:stemple}. 
See also~\cite{kn:bandelt,kn:stemple3,kn:zelinka}. 
These contain the Moore
graphs of diameter two, {\em i.e.\/}, the 5-cycle, 
the Petersen graph and 
the Hoffman-Singleton graph.%
\footnote{The only other possible 
degree that a Moore graph can have (thus with diameter 2 and girth 5) 
is 57. The existence of such a graph 
is unsettled~\cite{kn:miller}. The smallest unsettled case 
for a strongly regular graph with $\mu=1$ has parameters 
$(n,k,\lambda,\mu)=(400,21,2,1)$~\cite{kn:brouwer}.}  
Note that the Petersen graph has an 
induced $C_6$. 
It can be
shown (see~\cite{kn:kantor,kn:stemple,kn:scapellato,kn:blokhuis} 
and~\cite[Theorem~1.17.1]{kn:brouwer}) 
that, if $G$ is a geodetic
graph of diameter two then either 
\begin{enumerate}[\rm (i)]
\item $G$ contains a universal vertex,\footnote{In this case 
$G$ is a collection of cliques sharing one universal vertex.} or
\item $G$ is strongly regular, or
\item $G$ has exactly two vertex degrees $k_1>k_2$. Let $X_i$
be the set of vertices with degree $k_i$. Then $X_2$ is an
independent set. Every maximal clique that contains a vertex of
$X_1$ and $X_2$ has size two. Every maximal clique contained in
$X_1$ has size $k_1-k_2+2$. Furthermore, $n=k_1k_2+1$.
\end{enumerate}

Plesn\'ik~\cite{kn:plesnik} and Stemple~\cite{kn:stemple3} 
show that a geodetic graph $G$  
is homemorphic to $K_n$ if and only if 
there exists a function $f$ which assigns 
a nonnegative integer to every vertex of $K_n$ such that an edge 
$(x,y)$   
in $K_n$ has $f(x)+f(y)$ extra vertices 
in $G$. 
See~\cite[Section~7.3]{kn:buckley} for various constructions 
of geodetic graphs. 
A characterization of planar geodetic graphs appeared
in~\cite{kn:stemple2,kn:watkin}. 

\section{Acknowledgements}
%%%%%%%%%%%%%%%%%%%%%%%%%%

Ton~Kloks thanks his coauthors and the School of 
Computing of the University of Leeds for the pleasant collaboration.

\end{document}